\documentclass{article}

\usepackage{amsmath, amssymb}
\usepackage{graphicx, wrapfig} 
\usepackage{amsthm}
\usepackage{url}
\usepackage{color}
\usepackage{mathrsfs}

\theoremstyle{plain}
\newtheorem{theorem}{Theorem}

\theoremstyle{definition}
\newtheorem{definition}{Definition}
\newtheorem{example}{Example}

\newtheorem{proposition}{Proposition}
\newtheorem{lemma}{Lemma}
\newtheorem{corollary}{Corollary}

\newtheorem{remark}{Remark}

\title{Determinant of the crossing matrix of a braid}
\author{Ayaka Shimizu\thanks{Institute for Global Leadership, Ochanomizu University, Ohtsuka, Tokyo, 112-8610, Japan. Email: shimizu.ayaka@ocha.ac.jp, shimizu1984@gmail.com}}
\date{\today}

\begin{document}

\maketitle

\begin{abstract}
In this paper, we define a braid invariant, the purified determinant $P(b)$ of a braid $b$, considering the determinant of the crossing matrix of a pure braid derived from $b$, and show that $P(b_1 b_2)=P(b_2 b_1)$ for any pair of $n$-braids $b_1$ and $b_2$. 
\end{abstract}

\section{Introduction}

\begin{figure}[ht]
\centering
\includegraphics[width=9cm]{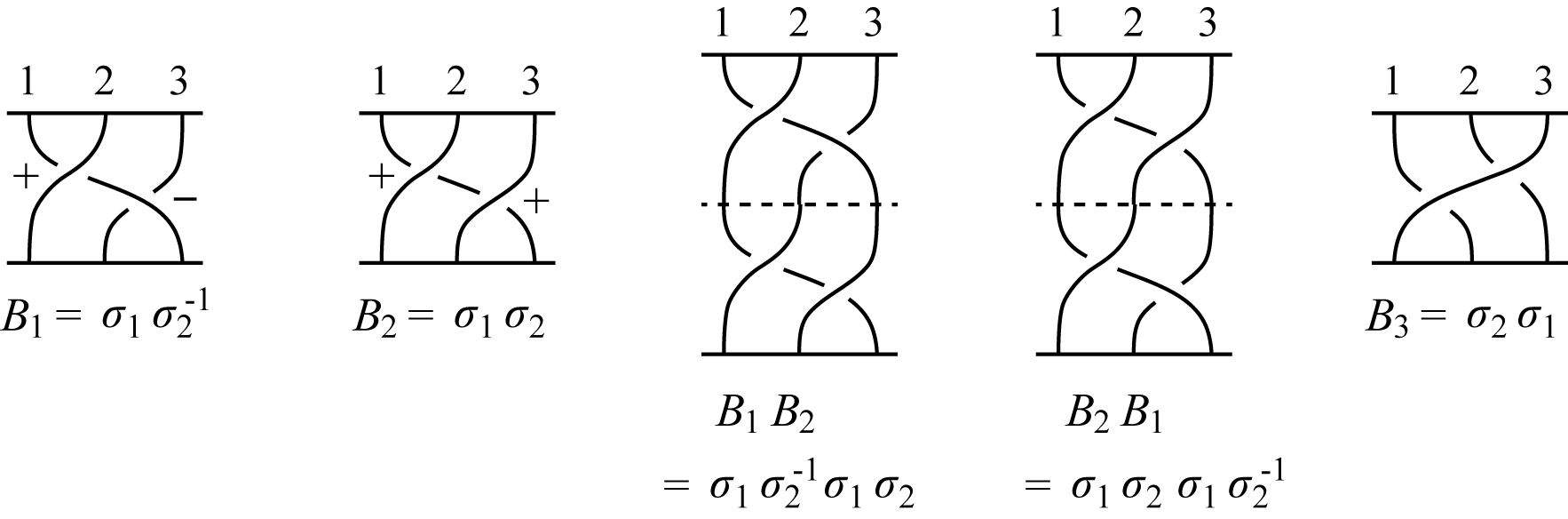}
\caption{Braid diagrams and their words. }
\label{f-braid}
\end{figure}

An {\it $n$-braid} (or a {\it braid}) $b$ is a set of $n$ disjoint strands in $\mathbb{R}^3$ such that the endpoints of each strand are attached to two horizontal bars and each strand runs from the upper bar to the lower bar monotonically. 
An {\it $n$-braid diagram} (or a {\it braid diagram}) $B$ is a regular projection of an $n$-braid on $\mathbb{R}^2$ as shown in Figure \ref{f-braid}. 
For each strand, we call the endpoint on the upper bar the {\it initial point} and the endpoint on the lower bar the {\it end point}.  
We call the strand whose initial point is at the $i^{th}$ position from the left the {\it $i^{th}$ strand}. \\

Each $n$-braid diagram is represented by a word that is generated by $\sigma_i^{\varepsilon}$, where $i \in \{ 1, 2, \dots , n-1 \}$ and $\varepsilon \in \{ +1, -1 \}$ as shown in Figure \ref{f-braid}. 
It is well known that two braid diagrams represent the same braid if and only if their words can be transformed into each other by a finite sequence of the following three types of equations: 
$\sigma_i \sigma_i^{-1} = \sigma_i^{-1} \sigma_i =1$, 
$\sigma_i \sigma_j= \sigma_j \sigma_i$ ($| j-i | >1$),  
$\sigma_i \sigma_j \sigma_i = \sigma_j \sigma_i \sigma_j$ ($| j-i | =1$). 
The last two equations are called {\it braid relations} (see \cite{Artin-1}, \cite{Artin-2}). 
For example, the two braid diagrams shown on the right-hand side in Figure \ref{f-braid} represent the same braid as $\sigma_1 \sigma_2 \sigma_1 \sigma_2^{-1} = \sigma_2 \sigma_1 \sigma_2 \sigma_2^{-1} =\sigma_2 \sigma_1$. 
We denote the set of all $n$-braids (respectively, $n$-braid diagrams) by $\mathfrak{B}_n$ (respectively, $\mathscr{B}_n$). 
For two $n$-braids $b_1$ and $b_2 \in \mathfrak{B}_n$ (respectively, $n$-braid diagrams $B_1$ and $B_2 \in \mathscr{B}_n$), the {\it product} $b_1 b_2$ (respectively, $B_1 B_2$) is defined to be the $n$-braid (respectively, $n$-braid diagram) that is obtained by connecting the end points of $b_1$ (respectively, $B_1$) to the initial points of $b_2$ (respectively, $B_2$) as shown in Figure \ref{f-braid}. 
We denote the product of $k$ braid diagrams $B$ by $B^k$. 

For a braid diagram $B= \sigma_{i_1}^{\varepsilon_1} \sigma_{i_2}^{\varepsilon_2} \dots \sigma_{i_k}^{\varepsilon_k} \in \mathscr{B}_n$ ($i_1, i_2, \dots , i_k \in \{ 1,2, \dots , n-1 \}$, $\varepsilon_1 , \varepsilon_2, \dots , \varepsilon_k \in \{ +1, -1 \}$), the {\it inverse} $B^{-1} \in \mathscr{B}_n$ is defined to be $B^{-1}= \sigma_{i_k}^{- \varepsilon_k} \dots \sigma_{i_2}^{- \varepsilon_2}  \sigma_{i_1}^{- \varepsilon_1} \in \mathscr{B}_n$. 
We note that $BB^{-1}$ represents the trivial $n$-braid. 
For a braid $b \in \mathfrak{B}_n$, the {\it inverse} $b^{-1} \in \mathfrak{B}_n$ is a braid such that $bb^{-1}$ is the trivial $n$-braid. 
For $B$ and $B' \in \mathscr{B}_n$, $B'$ is said to be {\it conjugate to} $B$ when $B'$ is expressed as $B'= A^{-1}BA$ for some $A \in \mathscr{B}_n$. 
Two braids $b$ and $b' \in \mathfrak{B}_n$ are said to be {\it conjugate in} $\mathfrak{B}_n$ when $b'$ is expressed as $b'=a^{-1}ba$ for some $a \in \mathfrak{B}_n$. 
We note that $b_1 b_2$ and $b_2 b_1$ are conjugate for any $b_1$, $b_2 \in \mathfrak{B}_n$ since $b_2 b_1 = b_1^{-1}(b_1 b_2)b_1$. 
A deformation from $b_1 b_2$ to $b_2 b_1$ or $B_1 B_2$ to $B_2 B_1$ is called a {\it conjugation}. \\

In this paper, we introduce a new braid invariant, the {\it purified determinant}, $P(b)$, of a braid $b \in \mathfrak{B}_n$ considering the determinant of the crossing matrix of a pure braid derived from $b$, and show that the value of $P$ is unchanged by a conjugation. 

\medskip
\begin{theorem}
If $b$ and $b' \in \mathfrak{B}_n$ are conjugate, then $P(b)=P(b')$. 
In particular, we have $P(b_1 b_2)=P(b_2 b_1)$ for any pair of $n$-braids $b_1$ and $b_2$. 
\label{thm-m1}
\end{theorem}
\medskip

The rest of the paper is organized as follows. 
In Section \ref{section-pre}, we review the definition of the crossing matrix and its properties.
In Section \ref{section-similar}, we discuss the similarity of the crossing matrices. 
In Section \ref{section-def}, we define the purified determinant for braid diagrams and braids. 
In Section \ref{section-proof}, we prove Theorem \ref{thm-m1}. 
In Section \ref{section-m2}, we discuss our future work.

\section{Crossing matrix}
\label{section-pre}

\begin{figure}[ht]
\centering
\includegraphics[width=5cm]{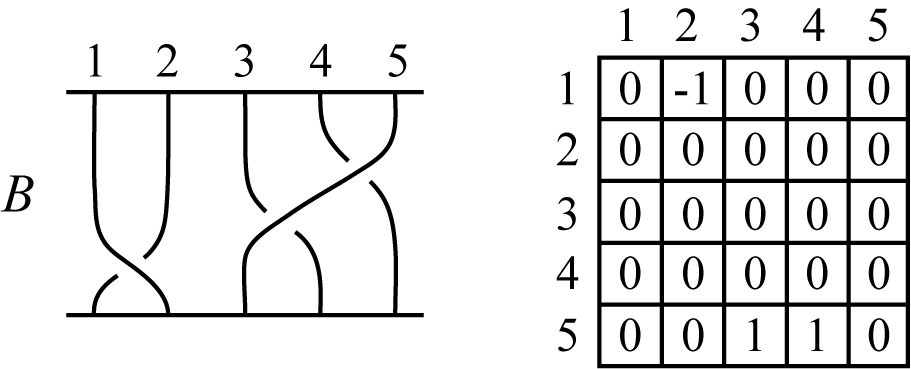}
\caption{A braid diagram $B$ and its crossing matrix $C(B)$. }
\label{f-rho}
\end{figure}

In this section, we review the crossing matrix and its properties. 
For a braid diagram, we call $\sigma_i$ a {\it positive crossing} and $\sigma_i^{-1}$ a {\it negative crossing}. 
For a matrix $M$, we denote the $(i,j)$-component by $M(i,j)$. 
The {\it crossing matrix} $C(B)=M$ of an $n$-braid diagram $B$ was defined by Burillo, Gutierrez, Krsti\'{c} and Nitecki in \cite{Bu} (see also \cite{Gu}) to be the $n \times n$ matrix such that $M(i,i)=0$ and $M(i,j)$ is the number of positive crossings minus negative crossings between the $i^{th}$ and $j^{th}$ strands where the $i^{th}$ strand is over the $j^{th}$ strand. (See Figure \ref{f-rho}.) 
The matrix $C(B)$ is unchanged under the braid relations, and the following proposition was shown in \cite{Bu}. 

\medskip 
\begin{proposition}[\cite{Bu}]
For each braid $b$, the crossing matrix $C(B)$ does not depend on the choice of the diagram $B$ of $b$. 
\label{prop-diag}
\end{proposition}
\medskip 

\noindent Then, the crossing matrix $C(b)$ of a braid $b$ is defined to be $C(b)=C(B)$ for any diagram $B$ of $b$. \\

For a braid $b \in \mathfrak{B}_n$ or a braid diagram $B \in \mathscr{B}_n$, when a strand has the initial point at the $i^{th}$ position from the left and has the end point at the $j^{th}$ position from the left, we denote $\rho (i)=j$ and $\rho (1,2, \dots , n )=(j_1 , j_2 , \dots , j_n )$ for $\rho (k )= j_k$. 
We call $\rho$ the {\it braid permutation} of $b$ or $B$. 
For example, the braid diagram $B$ in Figure \ref{f-rho} has the braid permutation $\rho$ as $\rho (1,2,3,4,5)=(2,1,4,5,3)$. 
We note that a braid and its diagram have the same braid permutation. 
By regarding $\rho$ as an element of the symmetric group, the identity $\rho^0 =id$, the inverse $\rho^{-1}$, the exponentiation $\rho^k$ for $k \in \mathbb{N}$ are also defined. 
For a permutation $\rho$, the {\it order}, $| \rho |$, of $\rho$ is defined to be the number $k \in \mathbb{N}$ such that $\rho^k =id$ and $\rho^l \neq id$ for any $l<k$. 
We say that a braid or a braid diagram is {\it pure} when the permutation is identity. 
The following proposition was shown in \cite{Bu}. 

\medskip
\begin{proposition}[\cite{Bu}]
An $n \times n$ matrix $M$ is the crossing matrix of some pure $n$-braid if and only if $M$ is a zero-diagonal symmetric integer matrix\footnote{
It has been conjectured that an $n \times n$ integer matrix $M$ is the crossing matrix of some {\it positive} pure $n$-braid if and only if $M$ is a non-negative integer zero-diagonal symmetric matrix with a property ``T0'' (\cite{Bu}), and addressed in \cite{Gu}, \cite{YAY}, \cite{AY-d}. 
}. 
\label{prop-pure}
\end{proposition}
\medskip

We define an operation on a square matrix that reorders the row and column simultaneously according to a permutation $\rho$. 

\medskip 
\begin{definition}
For an $n \times n$ matrix $M$ and a permutation $\rho$ on $(1,2, \dots ,n)$, we define $\rho (M)=N$ as $N(i,j)=M(\rho^{-1}(i), \rho^{-1}(j))$. 
\end{definition}
\medskip  

\begin{example}
For the crossing matrix $M$ and the braid permutation $\rho (1,2,3,4,5) = (2,1,4,5,3)$ of the braid diagram $B$ in Figure \ref{f-rho}, we have
\begin{align*}
M=
\begin{bmatrix}
0 & -1 & 0 & 0 & 0 \\
0 & 0 & 0 & 0 & 0 \\
0 & 0 & 0 & 0 & 0 \\
0 & 0 & 0 & 0 & 0 \\
0 & 0 & 1 & 1 & 0 
\end{bmatrix}
, \ \rho (M)=
\begin{bmatrix}
0 & 0 & 0 & 0 & 0 \\
-1 & 0 & 0 & 0 & 0 \\
0 & 0 & 0 & 1 & 1 \\
0 & 0 & 0 & 0 & 0 \\
0 & 0 & 0 & 0 & 0 
\end{bmatrix}
, \ \rho^2 (M)=
\begin{bmatrix}
0 & -1 & 0 & 0 & 0 \\
0 & 0 & 0 & 0 & 0 \\
0 & 0 & 0 & 0 & 0 \\
0 & 0 & 1 & 0 & 1 \\
0 & 0 & 0 & 0 & 0 
\end{bmatrix}
.
\end{align*}
\label{ex-rM}
\end{example}
\medskip

\noindent As mentioned in \cite{AY-5}, $M$ and $\rho (M)$ are similar, that is, there exists an invertible matrix $Q$ such that $\rho (M)= Q^{-1} MQ$, and therefore $\det (M) = \det (\rho (M))$. 
For two $n \times n$ matrices $M$ and $N$, we have $\rho (M+N)= \rho (M) + \rho (N)$ by definition. 
The following proposition was shown in \cite{Bu} for the crossing matrix of a braid product. 

\medskip 
\begin{proposition}[\cite{Bu}]
We have $C(B_1 B_2)=C(B_1)+ \rho (C(B_2))$, where $\rho$ is the braid permutation of $B_1$. 
\label{prop-C-product}
\end{proposition}
\medskip

\section{Similarity of crossing matrices}
\label{section-similar}

In this section, we discuss the similarity of crossing matrices of conjugate pure braids. 
We have the following lemma. 

\medskip 
\begin{lemma}
Let $B \in \mathscr{B}_n$, $i \in \{ 1, 2, \dots , n-1 \}$ and $\varepsilon \in \{ +1, -1 \}$. 
If $B$ is pure, then $C(B)$ and $ C(\sigma_i^{\varepsilon} B\sigma_i^{-\varepsilon})$ are similar to each other.
\label{lem-pure-det}
\end{lemma}
\medskip 

\noindent We remark that the order of the permutations of $B$ and $\sigma_i^{\varepsilon} B \sigma_i^{- \varepsilon}$ are the same since they are conjugate. 
Therefore, $B$ is pure if and only if $\sigma_i^{\varepsilon}B\sigma_i^{- \varepsilon}$ is pure (or consider the closure). \\

\begin{figure}[ht]
\centering
\includegraphics[width=4.5cm]{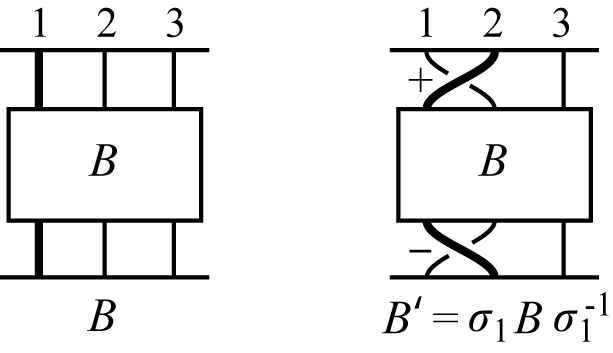}
\caption{Braid diagrams $B$ and $B'= \sigma_1 B \sigma_1^{-1}$. }
\label{f-conj}
\end{figure}

\medskip 
\noindent {\it Proof of Lemma \ref{lem-pure-det}}. 
Let $M_1 =C(B)$ and $M_2=C(B')$, where $B'= \sigma_i^{\varepsilon} B \sigma_i^{-\varepsilon}$. 
For $B'$, the $i^{th}$ (respectively, $(i+1)^{th}$) strand corresponds to the $(i+1)^{th}$ (respectively, $i^{th}$) strand of $B$ because there is a crossing $\sigma_i^{\varepsilon}$ of $B'$ between the $i^{th}$ and $(i+1)^{th}$ strands as the first crossing, which is above $B$. 
Also, the first crossing $\sigma_i^{\varepsilon}$ and the last crossing $\sigma_i^{- \varepsilon}$ of $B'$ are both between the $i^{th}$ and $(i+1)^{th}$ strands because $B'$ is pure, where the $(i+1)^{th}$ (respectively, $i^{th}$) strand is over when $\varepsilon =+1$ (respectively, $\varepsilon =-1$), and one of the two crossings is positive and the other is negative. 
(See Figure \ref{f-conj}.) 
Hence, we have $M_2 (i+1, i) = M_1 (i, i+1)$ and $M_2 (i,i+1)=M_1(i+1,i)$, and $M_2= \tau_{i, i+1} (M_1)$, where $\tau_{i, i+1}$ is a permutation that swaps the $i^{th}$ and $(i+1)^{th}$ elements. 
This implies that $M_1$ and $M_2$ are similar to each other. 
\hfill$\square$  \\
\medskip 

\noindent Since $C(B)$, $C( \sigma_{i_1}^{\varepsilon_1} B \sigma_{i_1}^{- \varepsilon_1})$, $C(\sigma_{i_2}^{\varepsilon_2} \sigma_{i_1}^{\varepsilon_1} B \sigma_{i_1}^{- \varepsilon_1} \sigma_{i_2}^{- \varepsilon_2} )$, $\dots$ are similar to each other when $B$ is pure by Lemma \ref{lem-pure-det}, we have the following corollary. 

\medskip 
\begin{corollary}
If a braid diagram $B' \in \mathscr{B}_n$ is conjugate to a pure braid diagram $B \in \mathscr{B}_n$, then $C(B)$ and $C(B')$ are similar to each other. 
\label{cor-pure-det}
\end{corollary}
\medskip

\noindent Lemma \ref{lem-pure-det} and Corollary \ref{cor-pure-det} do not hold for non-pure braid diagrams. 
For example, 
\begin{align*}
C(\sigma_1)= 
\begin{bmatrix}
0 & 0 & 0 \\
1 & 0 & 0 \\
0 & 0 & 0 
\end{bmatrix}
\text{ \ and } \ C(\sigma_2^{-1} \sigma_1 \sigma_2 )= 
\begin{bmatrix}
0 & 0 & 0 \\
1 & 0 & -1 \\
1 & 0 & 0 
\end{bmatrix}
\end{align*}
are not similar. 
Observe that the crossing matrices of them to the second power 
\begin{align*}
C((\sigma_1)^2)= 
\begin{bmatrix}
0 & 1 & 0 \\
1 & 0 & 0 \\
0 & 0 & 0 
\end{bmatrix}
\text{ \ and } \ C((\sigma_2^{-1} \sigma_1 \sigma_2 )^2)= 
\begin{bmatrix}
0 & 0 & 1 \\
0 & 0 & 0 \\
1 & 0 & 0 
\end{bmatrix}
\end{align*}
are similar to each other. 
We have the following lemma for non-pure braid diagrams. 

\medskip
\begin{lemma}
If $B'\in \mathscr{B}_n$ is conjugate to $B \in \mathscr{B}_n$, then $C(B^r)$ and $C((B')^r)$ are similar to each other, where $r$ is the order of the braid permutation of $B$. 
\label{lem-P}
\end{lemma}
\medskip 

\begin{proof}
Let $B'=A^{-1}BA$. 
Since $B^r$ is pure and $(B')^r=(A^{-1}BA) (A^{-1}BA) \dots (A^{-1}BA)$ has the same crossing matrix to $A^{-1} B^r A$ by Proposition \ref{prop-diag}, $B^r$ and $(B')^r$ have similar crossing matrices by Corollary \ref{cor-pure-det}. 
\end{proof}
\medskip 

\noindent We note that the order of the permutation of $B'$ is also $r$ since $B'$ is conjugate to $B$. 
For braids, we have the following proposition. 

\medskip
\begin{proposition}
If $b$ and $b' \in \mathfrak{B}_n$ are conjugate, then $C(b^r)$ and $C((b')^r)$ are similar to each other, where $r$ is the order of the permutation of $b$. 
\label{prop-similar}
\end{proposition}
\medskip 

\begin{proof}
Take braid diagrams of $b$ and $b'$ in the form of $B$ and $B'= A^{-1}BA \in \mathscr{B}_n$. 
By Proposition \ref{prop-diag} and Lemma \ref{lem-P},  the crossing matrices $C(b^r)=C(B^r)$ and $C((b')^r) = C((B')^r)$ are similar to each other. 
\end{proof}

\section{Definition of the purified determinant}
\label{section-def}

In this section, we define the purified determinant and see examples. 

\medskip 
\begin{definition}
For an $n$-braid diagram $B \in \mathscr{B}_n$, let $\rho$ be the braid permutation of $B$. 
The {\it purified determinant} of $B$, denoted by $P(B)$, is defined to be $P(B)= \det (C(B^{| \rho |}))$. 
\end{definition}
\medskip 

\noindent We note that $B^{| \rho |}$ is a pure braid diagram that is obtained by taking the product of  $| \rho |$ braid diagrams $B$. 
The number of crossings of $B^{| \rho |}$ may be large. 
The following proposition implies that we can also calculate the value of $P(B)$ by considering the crossing matrix of $B$ instead of $B^{| \rho |}$. 

\medskip 
\begin{proposition}
\begin{align*}
P(B)= \det \left( \sum_{k=0}^{| \rho | -1} \rho^k (C(B)) \right) .
\end{align*}
\end{proposition}
\medskip  

\begin{proof}
By Proposition \ref{prop-C-product}, we have 
\begin{align*}
C(B^{| \rho |}) & = C(B)+ \rho (C(B^{| \rho |-1})) \\
& = C(B)+ \rho ( C(B) + \rho (C(B^{| \rho |-2}))) \\
& = C(B) + \rho (C(B))+ \rho^2 (C(B^{| \rho |-2})) \\
& = C(B) + \rho (C(B))+ \rho^2 (C(B)) + \dots + \rho^l (C(B^{| \rho |-l})).
\end{align*}

\end{proof}
\medskip 

\begin{example}
For the 5-braid diagram $B$ in Figure \ref{f-rho}, we have $| \rho | =6$, 
\begin{align*}
\sum_{k=0}^5 \rho^k \left(
\begin{bmatrix}
0 & -1 & 0 & 0 & 0 \\
0 & 0 & 0 & 0 & 0 \\
0 & 0 & 0 & 0 & 0 \\
0 & 0 & 0 & 0 & 0 \\
0 & 0 & 1 & 1 & 0 
\end{bmatrix}
\right) = 
\begin{bmatrix}
0 & -3 & 0 & 0 & 0 \\
-3 & 0 & 0 & 0 & 0 \\
0 & 0 & 0 & 2 & 2 \\
0 & 0 & 2 & 0 & 2 \\
0 & 0 & 2 & 2 & 0 
\end{bmatrix}
\end{align*}
and $P(B)=-144$. 
(See Example \ref{ex-rM} for $\rho (C(B))$ and $\rho^2 (C(B))$.)
\end{example}
\medskip

\begin{remark}
Since the crossing matrix $C(B)$ is a braid invariant $C(b)$, the purified determinant is also a braid invariant. 
Namely, we can define the purified determinant of $b \in \mathfrak{B}_n$ as $P(b)=P(B)$ for any diagram $B$ of $b$. 
\label{rem-inv}
\end{remark}
\medskip

\begin{example}
A {\it fundamental braid}, $\Delta_n$, is a positive $n$-braid that is obtained from a trivial braid by turning over the lower bar. 
We have 
$$P( \Delta_n) = - (-1)^n (n-1).$$
\end{example}
\medskip 

\begin{proof}
The order of permutation of each $\Delta_n$ is 2, and each $\Delta_n$ has an $n$-braid diagram such that each pair of strands has exactly one positive crossing. 
Hecne, the matrix $M_n=C((\Delta_n)^2)$ is a matrix such that $M_n(i,i)=0$ and $M_n(i,j)=1$ for $i \neq j$. 
To calculate $\det (M_n)$, subtract the $2^{nd}$ column from the $1^{st}$ column, and apply a Laplace expansion along the $1^{st}$ column. 
Then, break down the $1^{st}$ row of the $2^{nd}$ term, and we obtain a recurrence formula $\det (M_n)=-2 \det (M_{n-1})- \det (M_{n-2})$. 
For example, when $n=5$, we have 
\begin{align*}
&
\begin{vmatrix}
0 & 1 & 1 & 1 & 1 \\
1 & 0 & 1 & 1 & 1 \\
1 & 1 & 0 & 1 & 1 \\
1 & 1 & 1 & 0 & 1 \\
1 & 1 & 1 & 1 & 0 
\end{vmatrix}
=
\begin{vmatrix}
-1 & 1 & 1 & 1 & 1 \\
1 & 0 & 1 & 1 & 1 \\
0 & 1 & 0 & 1 & 1 \\
0 & 1 & 1 & 0 & 1 \\
0 & 1 & 1 & 1 & 0 
\end{vmatrix}
= -
\begin{vmatrix}
0 & 1 & 1 & 1 \\
1 & 0 & 1 & 1 \\
1 & 1 & 0 & 1 \\
1 & 1 & 1 & 0 
\end{vmatrix}
-
\begin{vmatrix}
1 & 1 & 1 & 1 \\
1 & 0 & 1 & 1 \\
1 & 1 & 0 & 1 \\
1 & 1 & 1 & 0 
\end{vmatrix}\\
& =-
\begin{vmatrix}
0 & 1 & 1 & 1 \\
1 & 0 & 1 & 1 \\
1 & 1 & 0 & 1 \\
1 & 1 & 1 & 0 
\end{vmatrix}
- \left( 
\begin{vmatrix}
0 & 1 & 1 & 1 \\
1 & 0 & 1 & 1 \\
1 & 1 & 0 & 1 \\
1 & 1 & 1 & 0 
\end{vmatrix}
+
\begin{vmatrix}
1 & 0 & 0 & 0 \\
1 & 0 & 1 & 1 \\
1 & 1 & 0 & 1 \\
1 & 1 & 1 & 0 
\end{vmatrix}
\right) 
= -2 \begin{vmatrix}
0 & 1 & 1 & 1 \\
1 & 0 & 1 & 1 \\
1 & 1 & 0 & 1 \\
1 & 1 & 1 & 0 
\end{vmatrix}
-
\begin{vmatrix}
0 & 1 & 1 \\
1 & 0 & 1 \\
1 & 1 & 0 
\end{vmatrix}.
\end{align*}
\end{proof}
\medskip

\section{Proof of the main theorem and examples}
\label{section-proof}

In this section, we prove Theorem \ref{thm-m1}. \\

\medskip
\noindent {\it Proof of Theorem \ref{thm-m1}.} 
If $b$ and $b' \in \mathfrak{B}_n$ are conjugate, then the crossing matrices $C(b^r)$ and $C((b')^r)$ are similar to each other by Proposition \ref{prop-similar}, where $r$ is the order of the braid permutation of $b$ and $b'$. 
Hence, their determinants are equal, namely, $P(b)=P(b')$. 
\hfill$\square$  \\
\medskip 

\noindent In the same way, we have the following corollary from Proposition \ref{prop-similar}. 

\medskip 
\begin{corollary}
The rank, eigenvalues, and characteristic polynomial of $C(b^{ | \rho | })$ are also braid invariants that are unchanged under conjugation. 
\end{corollary}
\medskip

\noindent We note that the matrix $C(b^{| \rho |})$ definitely has real eigenvalues for any braid $b$ since the crossing matrix of the pure braid $b^{| \rho |}$ is symmetric by Proposition \ref{prop-pure}. 

\medskip 
\begin{example}
The braid diagrams $B$ and $B' \in \mathscr{B}_3$ in Figure \ref{f-P-ex} represent different braids since they have different crossing matrices. 
Moreover, they are not conjugate since they have the different values of purified determinant\footnote{Eigenvalues of $C(B^2)$ are $1$,  $\frac{-1 \pm \sqrt{33}}{2}$ and that of $C((B')^2)$ are $-1$, $\frac{1 \pm \sqrt{33}}{2}$. This also can conclude that $B$ and $B'$ are not conjugate. } as $P(B)=-8$ and $P(B')=8$, even the order of their permutations are the same and the linking number of the closures (see Section \ref{section-m2}) are the same. 
\begin{figure}[ht]
\centering
\includegraphics[width=5.5cm]{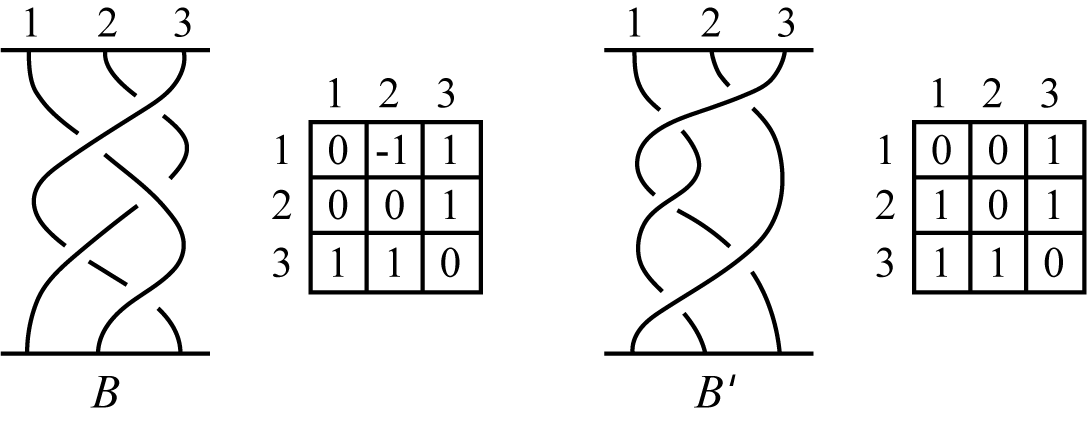}
\caption{Braid diagrams $B$ and $B'$ with the purified determinants $P(B)=-8$ and $P(B')=8$. }
\label{f-P-ex}
\end{figure}
\end{example}
\medskip

A braid diagram $B$ is said to be {\it layered} when the set of all the strands of $B$ can be divided into two sets $S_1$ and $S_2$ so that any strand in $S_1$ does not pass under a strand in $S_2$. 
We denote a braid diagram that is obtaind from $B$ by ignoring all the strands in $S_2$ (respectively, $S_1$) by $B_1$ (respectively, $B_2$) and denote as $B=B_1 \oplus B_2$. 
For pure braid diagrams, we have the following proposition in the same way to Theorem 1 in \cite{AY-5}. 

\medskip 
\begin{proposition}
When a layered braid diagram $B=B_1 \oplus B_2$ is pure, we have $P(B)=P(B_1)P(B_2)$. 
\end{proposition}
\medskip 

\begin{proof}
Take a sequence $\mathbf{s}=( s^1_1, s^1_2, \dots , s^1_{k_1}, s^2_1, s^2_2, \dots , s^2_{k_2})$ of the strands of $B$ so that $s^1_1, s^1_2, \dots , s^1_{k_1} \in S_1$ and $s^2_1, s^2_2, \dots , s^2_{k_2} \in S_2$. 
Reorder the row and column of $C(B)$ simultaneously according to the order of $\mathbf{s}$ to obtain a matrix $M$, which is similar to $C(B)$. 
Observe that $M$ is in the form of 
\begin{align*}
M=
\begin{bmatrix}
M_1 & O \\
O & M_2
\end{bmatrix}
\end{align*}
with the blocks $M_1 =C(B_1)$ and $M_2 = C(B_2)$. 
We note that $C(B)$ is symmetric by Proposition \ref{prop-pure} and $M$ is also symmetric. 
The block on the lower left is $O$ because no strand in $S_2$ passes over a strand in $S_1$. 
The block on the upper right is also $O$ by symmetry. \\
Since $B$ is pure, we have $P(B)= \det (C(B))$, and $\det (C(B))= \det (M) = \det (M_1) \det (M_2) = P(B_1) P(B_2)$. 
We note that $B_1$ and $B_2$ are also pure, and therefore $P(B_i)= \det (C(B_i))= \det (M_i)$. 
\end{proof}
\medskip

By definition, a braid diagram $B$ with order of the permutation $r$ and the pure braid diagram $B^r$ have the same value of the purified determinant. 
We define a more detailed invariant to distinguish such $B$ and $B^r$.  

\medskip 
\begin{definition}
For a braid diagram $B$ and its permutation $\rho$, let $PP(B)=(| \rho | , P(B))$. 
We call $PP(B)$ the {\it P-pair} of $B$. 
\end{definition}
\medskip 

\noindent By Theorem \ref{thm-m1}, we have the following proposition. 

\medskip 
\begin{proposition}
The P-pair is a braid invariant. 
Thai is, we can define the P pair of a braid $b \in \mathfrak{B}_n$ as $PP(b)=PP(B)$ for any diagram $B$ of $b$. 
Moreover, $PP(b)$ is unchanged under conjugation. 
\end{proposition}
\medskip 

\begin{example}
Let $B \in \mathscr{B}_3$ be a braid diagram represented by $\sigma_1^{-1} \sigma_2$. 
Although $B$ and $B^3$ have the same value of the purified determinant by definition, they have different P-pairs as $PP(B)=(3,0)$ and $PP(B^3)=(1,0)$. 
\end{example}

\section{Future work}
\label{section-m2}

In this section, we discuss braid invariants derived from the crossing matrix with our future goal of constructing a link invariant based on the crossing matrix. 
From each braid diagram $B$, we obtain an oriented link diagram by connecting each pair of the $i^{th}$ end point and the $i^{th}$ initial point by an arc without creating crossings. 
We call it the {\it closure of} $B$ and denote it by $\widehat{B}$. 
It was shown in \cite{Alex} by Alexander that any oriented link can be represented by a braid diagram $B \in \mathscr{B}_n$ for some $n$ as the closure $\widehat{B}$. \\

\begin{figure}[ht]
\centering
\includegraphics[width=7cm]{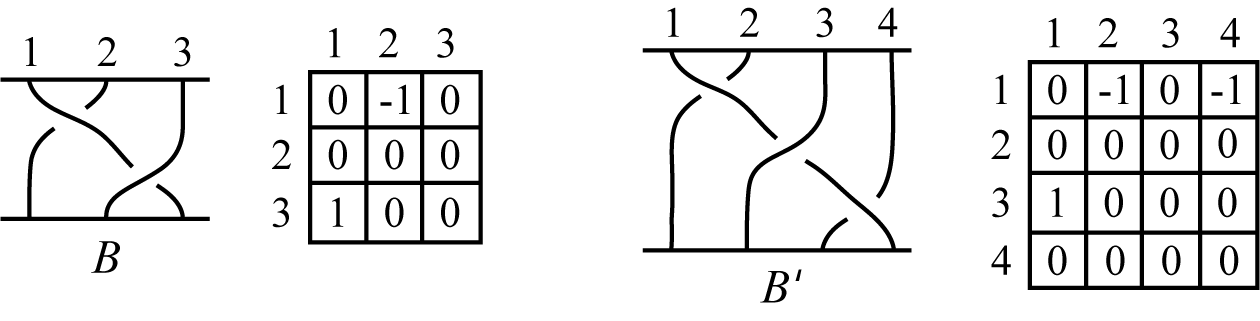}
\caption{The braid diagram $B'$ is obtained from $B$ by a stabilization. }
\label{f-st}
\end{figure}

For an $n$-braid diagram $B \in \mathscr{B}_n$, a  {\it stabilization} is the deformation to obtain $B \sigma_n^{\varepsilon} \in \mathscr{B}_{n+1}$, where $\varepsilon \in \{ +1, -1 \}$. 
The opposite deformation is called a {\it destabilization}. 
Markov showed in \cite{Markov} that for two braids $b$ and $b'$, the closures $\widehat{b}$ and $\widehat{b'}$ represent the same link type if and only if $b$ and $b'$ are related by a finite number of conjugations, stabilizations and destabilizations. 
(See, for example, \cite{BM}, \cite{K-s}.) 
We note that the purified determinant can be changed by a stabilization. 
For example, we have $P(B)=0$ and $P(B')=-3$ for the braid diagrams in Figure \ref{f-st}. 
Now we define a braid invariant that is unchanged under a stabilization. 

\medskip 
\begin{definition}
For a braid diagram $B \in \mathscr{B}_n$, let $Q(B)= \det ( C(B) +I)$, where $I$ is the $n \times n$ identity matrix. 
We call $Q(B)$ the {\it Q-determinant} of $B$. 
\end{definition}
\medskip 

\noindent For example, we have $Q(B)=1$, $Q(B')=1$ for the braid diagrams in Figure \ref{f-st}. 

\medskip 
\begin{remark}
Since the crossing matrix is a braid invariant, the Q-determinant is also a braid invariant. 
Namely, the Q-determinant of a braid $b$ is defined to be $Q(b)=Q(B)$ for any diagram $B$ of $b$. 
\end{remark}
\medskip 

\begin{proposition}
Let $B' \in \mathscr{B}_{n+1}$ be a braid diagram which is obtained from a braid diagram $B \in \mathscr{B}_n$ by a stabilization, namely $B'=B \sigma_n^{\varepsilon}$, where $\varepsilon \in \{ +1, -1 \}$. 
Then $Q(B')=Q(B)$. 
\end{proposition}
\medskip 

\begin{proof}
Let $M=C(B)$, $M'=C(B')$. 
Let $I$ (respectively, $I'$) be the $n \times n$ (respectively, $(n+1) \times (n+1)$) identity matrix. 
Suppose that the crossing $\sigma_n^{\varepsilon}$ of $B'$ at the lower right is a crossing between the $k^{th}$ and $(n+1)^{th}$ strands. 
\begin{itemize}
\item[(1)] Suppose that $\varepsilon =+1$. 
We have \\
$\bullet$ $M'(i,j)=M(i,j)$ when $i, j \leq n$, \\
$\bullet$ $M'(n+1, k)=1$, \\
$\bullet$ $M'(n+1, j)=0$ when $j \neq k$, and \\
$\bullet$ $M'(i, n+1)=0$ for any $i$. \\
Apply $M'+I'$ the Laplace expansion along the $(n+1)^{th}$ column to obtain $\det (M' +I') = \det (M +I)$. 
\item[(2)] Suppose that $\varepsilon =-1$. 
We have \\
$\bullet$ $M'(i,j)=M(i,j)$ when $i, j \leq n$, \\
$\bullet$ $M'(k, n+1)=-1$, \\
$\bullet$ $M'(i, n+1)=0$ when $i \neq k$, and \\
$\bullet$ $M'(n+1, j)=0$ for any $j$. \\
Apply $M'+I'$ the Laplace expansion along the $(n+1)^{th}$ row to obtain $\det (M' +I') = \det (M +I)$. 
\end{itemize}
Hence, $Q(B')=Q(B)$. 
\end{proof}
\medskip 

\begin{corollary}
The Q-determinant of a braid is unchanged under stabilizations and destabilizations. 
\end{corollary}
\medskip 

\noindent We note that the Q-determinant can be changed under conjugations. 
For example, we have $Q( \sigma_2^{-3} \sigma_1)=-1$ and $Q( \sigma_2^{-2} \sigma_1 \sigma_2^{-1})=1$. 
Employing the crossing matrix, constructing a new braid invariant that is unchanged under both conjugations and stabilizations would be our future goal. 
It would also be a link invariant which distinguishes oriented or unoriented links from a new aspect.

\section*{Acknowledgment}
This work was partially supported by the JSPS KAKENHI Grant Number JP21K03263. 
The author thanks all the members and participants in the Maebashi Braid Seminar for valuable discussions on braids. 
In particular, she is very grateful to Yoshiro Yaguchi for helpful suggestions and valuable discussions.

\end{document}